\documentclass[12pt]{article}
\usepackage{amssymb,amsmath,amsfonts,amsthm,amscd,latexsym,verbatim,graphics,epsfig,indentfirst}
\usepackage{geometry}
\def\id{\mathrm{id}}
\def\dim{\mathrm{dim}\,}
\def\pre{\mathrm{pre}}
\def\post{\mathrm{post}}
\def\RB{\mathrm{RB}}
\def\Var{\mathrm{Var}}
\def\As{\mathrm{As}}
\def\Lie{\mathrm{Lie}}

\def\lbar{\prec }
\def\rbar{\succ }

\begin{document}

\begin{flushright}
""Š 512.579
\end{flushright}

\begin{center}
{\Large
Universal enveloping Rota---Baxter algebras of
preassociative and postassociative algebras}

V.Yu. Gubarev
\end{center}

\begin{abstract}
Universal enveloping Rota---Baxter algebras of
preassociative and postassociative algebras
are constructed. The question of Li Guo is answered:
the pair of varieties $(\RB_\lambda\As,\post\As)$
is a PBW-pair and the pair $(\RB\As,\pre\As)$ is not.

\end{abstract}

\section*{Introduction}

Linear operator $R$ defined on an algebra $A$ over the key field $\Bbbk$
is called Rota---Baxter operator (RB-operator, for short) of a weight $\lambda\in\Bbbk$
if it satisfies the relation
\begin{equation}\label{eq:RB}
R(x)R(y) = R( R(x)y + xR(y) + \lambda xy), \quad x,y\in A.
\end{equation}

An algebra with given RB-operator acting on it
is called Rota---Baxter algebra (RB-algebra, for short).

G. Baxter defined (commutative) RB-algebra in 1960 \cite{Baxter60},
solving an analytic problem. The relation \eqref{eq:RB} with $\lambda = 0$
appeared as a generalization of by part integration formula.
J.-C. Rota and others \cite{Rota68,Cartier72} studied combinatorial properties
of RB-operators and RB-algberas.
In 1980s, the deep connection between Lie RB-algebras and Young---Baxter equation
was found \cite{BelaDrin82,Semenov83}.
To the moment, there are a lot of applications of RB-operators in mathematical physics,
combinatorics, number theory, and operad theory \cite{QFT00,QFT01,NumThe,GubKol2013}.

There exist different constructions of free commutative RB-algebra,
see the articles of J.-C. Rota, P. Cartier, and L. Guo \cite{Rota68,Cartier72,GuoKeigher}.
In 2008, K. Ebrahimi-Fard and L. Guo obtained
free associative RB-algebra \cite{FardGuo08}.
In 2010, L.A. Bokut et al \cite{BokutChen} got a linear base of free associative RB-algebra
with the help of Gr\"{o}bner---Shirshov technique.
Diverse linear bases of free Lie RB-algebra
were recently found in \cite{Gub2016,GubKol2016,Chen16}.

Pre-Lie algebras were introduced in 1960s independently
by E.B. Vinberg, M. Gerstenhaber, and J.-L. Koszul \cite{Vinberg63,Gerst63,Koszul61},
they satisfy the identity $(x_1 x_2)x_3 - x_1(x_2 x_3) = (x_2 x_1) x_3 - x_2 (x_1 x_3)$.

J.-L. Loday \cite{Dialg99} defined the notion of (associative) dendriform dialgebra,
we will call it preassociative algebra or associative prealgebra.
Preassociative algebra is a vector space with two bilinear operations $\succ,\prec$
satisfying the identities
$$
\begin{gathered}
(x_1\succ x_2+x_1\prec x_2)\succ x_3 = x_1\succ (x_2 \succ x_3), \quad
(x_1\succ x_2)\prec x_3=x_1\succ(x_2\prec x_3), \\
x_1\prec(x_2\succ x_3+x_2\prec x_3)=(x_1\prec x_2)\prec x_3.
\end{gathered}
$$

In \cite{Loday95}, J.-L. Loday also defined Zinbiel algebra
(we will call it as precommutative algebra),
on which the identity
$(x_1\succ x_2 + x_2\succ x_1)\succ x_3 = x_1\succ (x_2\succ x_3)$
holds. Any preassociative algebra
with the identity $x\succ y = y\prec x$ is a precommutative algebra
and with respect to the new operation $x\cdot y = x\succ y - y\prec x$
is a pre-Lie algebra.
There is an open problem if every pre-Lie algebra
injectively embeds into its universal enveloping preassociative algebra
(in affirmative case, will be the pair of varieties of pre-Lie
and preassociative algebras a PBW-pair \cite{PBW}?).

In \cite{Trialg01}, there was also defined
(associative) dendriform trialgebra, i.e.,
an algebra with the operations $\prec,\succ,\cdot$ satisfying
certain 7 axioms.
(We will call such algebra as postassociative algebra or associative postalgebra.)
Post-Lie algebra \cite{Vallette2007} is an algebra
with two bilinear operations $[,]$ and $\cdot$; moreover,
Lie identities with respect to $[,]$ hold
and the next identities are satisfied:
$$
(x \cdot y) \cdot z - x \cdot (y \cdot z)
- (y \cdot x) \cdot z + y \cdot (x \cdot z) = [y,x]\cdot z, \quad
x\cdot [y,z] = [x \cdot y,z] + [y,x\cdot z].
$$

Given a binary operad $\mathcal P$,
the notion of succesor \cite{BBGN2012} provides
the defining identities for pre- and post-$\mathcal P$-algebras.
Equivalently, one can define the operad of pre- and post-$\mathcal P$-algebras
as $\mathcal P \bullet \mathrm{Pre}\Lie$ and
$\mathcal P \bullet \mathrm{Post}\Lie$ respectively.
Here $\mathrm{Pre}\Lie$ denotes the operad of pre-Lie algebras
and $\mathrm{Post}\Lie$ --- the operad of post-Lie algebras,
$\mathcal{V}\bullet\mathcal{W}$ is the black product of operads
$\mathcal{V}$, $\mathcal{W}$
(see \cite{GinzKapr} about operads and Manin products).
Hereinafter, pre- and postalgebra will denote
pre- and post-$\mathcal P$-algebra for some variety (operad) $\mathcal P$.

In 2000, M. Aguiar \cite{Aguiar00} remarked that
any associative algebra with defined on it a Rota---Baxter operator $R$
of zero weight with respect to the operations
$a\lbar b=aR(b)$, $a \rbar b=R(a)b$
is a preassociative algebra. In 2002, K. Ebrahimi-Fard \cite{Fard02}
showed that one can additionally define on an associative RB-algebra
of nonzero weight $\lambda$ the third operation
$a\cdot b = \lambda ab$ and get the structure of postassociative algebra
under the operations $\lbar$, $\rbar$, $\cdot$.
In 2007, the notion of universal enveloping RB-algebras
of pre- and postassociative algebra was introduced \cite{FardGuo07}.

For free preassociative algebra $C$,
injectivity of embedding $C$ into its
universal enveloping was proved in \cite{FardGuo07}.

In 2010, with the help of Gr\"{o}bner---Shirshov bases \cite{BokutChen},
Yu. Chen and Q. Mo proved that any preassociative algebra
over the field of zero characteristic injectivily emdebs into
its universal enveloping RB-algebra \cite{Chen11}.

In 2010, C. Bai et al \cite{BaiGuoNi10} introduced $\mathcal{O}$-operator,
a generalization of RB-operator, and stated that
any associative pre- and postalgebra injectivily embeds
into an algebra with $\mathcal{O}$-operator.

In \cite{BBGN2012}, the construction of M. Aguiar and K. Ebrahimi-Fard
was generalized on the case of arbitrary operad, not only associative.

In \cite{GubKol2013}, given a variety $\Var$,
it was proved that any pre-$\Var $-algebra
injectivily embeds into its universal enveloping
$\Var $-RB-algebra of zero weight and any post-$\Var $-algebra
injectivily embeds into its universal enveloping
$\Var $-RB-algebra of weight $\lambda\neq0$.
Based on the last results, we have

\textbf{Problem 1}.
To construct an universal enveloping RB-algebra of pre- and postalgebra.

In the comments to the head V of the unique monograph on RB-algebras \cite{Guo2011},
L. Guo actually stated the following

\textbf{Problem 2}.
To clarify if the pairs of varieties $(\RB\As,\mathrm{pre}\As)$
and $(\RB_\lambda\As,\mathrm{post}\As)$ for $\lambda\neq0$
are PBW-pairs \cite{PBW}.

Here $\RB\As$ ($\RB_\lambda\As$) denotes the variety of associative algebras
endowed with an RB-operator of (non)zero weight $\lambda$.

In the current article, Problem 1 is solved for associative pre- and postalgebras
and Problem 2 is solved completely.

In \S1, we state preliminaries about RB-algebras, PBW-pairs,
preassociative and post\-asso\-ci\-ative algebras.
Universal enveloping RB-algebras
of preassociative (\S2) and post\-asso\-ciative algebras
(\S3) are constructed.
As a corollary, we obtain that the pair of varieties
$(\RB_\lambda\As,\post\As)$
is a PBW-pair and the pair $(\RB\As,\pre\As)$ is not.

\section{Preliminaries}

\subsection{RB-operator}
Let us consider some well-known examples of RB-operators (see, e.g., \cite{Guo2011}):

{\bf Example 1}.
Given an algebra $A$ of continious functions on $\mathbb{R}$,
an integration operator
$R(f)(x) = \int\limits_{0}^x f(t)\,dt$
is an RB-operator on $A$ of zero weight.

{\bf Example 2}.
Given an invertible derivation $d$ on an algebra $A$,
$d^{-1}$ is an RB-operator on $A$ of zero weight.

{\bf Example 3}.
Let $A = A_1 \oplus A_2$ be a direct sum (as vector space)
of its two subalgebras. An operator $R$ defined as
$$
R|{}_{A_1}\equiv 0,\quad R|{}_{A_2}\equiv \lambda\id,
$$
where $\id$ denotes the identical map,
is an RB-operator on $A$ of the weight $-\lambda$.

Further, unless otherwise specified,
RB-operator will mean RB-operator of zero weight.
We denote a free algebra of a variety $\Var$
generated by a set $X$ by
$\Var\langle X\rangle$, and a free RB-algebra
of a weight $\lambda$ respectively by
$\RB_\lambda\Var\langle X\rangle$.
For short, denote
$\RB_0\Var\langle X\rangle$ by
$\RB\Var\langle X\rangle$.

\subsection{PBW-pair of varities}

In 2014, A.A. Mikhalev and I.P. Shestakov
introduced the notion of a PBW-pair \cite{PBW};
it generalizes the relation between
the varieties of associative and Lie algebras
in the spirit of Poincar\'{e}---Birkhoff---Witt theorem.

Given varieties of algebras
$\mathcal{V}$ and $\mathcal{W}$, let
$\psi\colon \mathcal{V}\to \mathcal{W}$
be a such functor that maps an algebra $A\in \mathcal{V}$
to the algebra $\psi(A)\in \mathcal{W}$
preserving $A$ as vector space but changing the operations on $A$.
There exists left adjoint functor to the $\psi$
called universal enveloping algebra and denoted as $U(A)$.
Defining on $U(A)$ a natural ascending filtration,
we get associated graded algebra $\mathrm{gr}\,U(A)$.

A pair of varieties $(\mathcal{V},\mathcal{W})$
with the functor $\psi\colon \mathcal{V}\to \mathcal{W}$
is called PBW-pair if $\mathrm{gr}\,U(A)\cong U(\mathrm{Ab}\,A)$.
Here $\mathrm{Ab}\,A$ denotes the vector space $A$
with trivial multiplicative operations.

\subsection{Free associative RB-algebra}

In \cite{FardGuo08}, free associative RB-algebra
in the terms of rooted forests and trees was constructed.
In \cite{BokutChen} Gr\"{o}bner---Shirshov theory
of associative RB-algebras was developed.
Based on the results \cite{BokutChen,FardGuo08},
one can get by induction linear base of free associative RB-algebra
$\RB_\lambda\As\langle X\rangle$ generated by a set $X$,
where $\lambda$ is either zero or nonzero. The base consists of
{\it RB-associative words} of the form
$$
w_1 R(u_1)w_2 R(u_2)\ldots w_k R(u_k)w_{k+1},
$$
where $w_i\in S(X)$, $i=2,\ldots,k$,
are elements of free semigroup generated by the set $X$,
$w_1,w_{k+1}\in S(X)\cup\emptyset$,
$u_1,\ldots,u_k$ are RB-associative words obtained on the previous inductive step.
The product of such words $u\cdot v$ differs from the concatenation $uv$
only if $u,v$ have the following view: $u = u_1 R(u_2)$, $v = R(v_1) v_2$.
In this case,
$$
u\cdot v = u_1 R(R(u_2)\cdot v_1 + u_2\cdot R(v_1) + \lambda u_2\cdot v_1)v_2.
$$
Further, we will refer to the constructed base as the (standard) base of $\RB_\lambda\As\langle X\rangle$.

Given a word $u$ from the standard base of $\RB_\lambda\As\langle X\rangle$,
the number of appearances of the symbol $R$ in the notation of $u$
is called $R$-degree of the word $u$, denotation: $\deg_R(u)$.

We will call any element of the standard base of $\RB_\lambda\As\langle X\rangle$
of the form $R(w)$ as $R$-letter.
By $X_\infty$ we denote the union of the set $X$ and the set of all $R$-letters.
Let us define a degree $\deg u$ of the word $u$ from the standard base
as the length of $u$ in the alphabet $X_\infty$.
For example, given the word $u = x_1 R(x_2)x_3 R(x_4 x_5)$,
we have $\deg_R(u) = 2$ and $\deg u = 4$, since $u$
consists of four different letters $x_1$, $x_3$, $R(x_2)$, $R(x_4 x_5)\in X_\infty$.

\subsection{Preassociative algebra}

A linear space within two bilinear operations $\succ,\prec$
is called preassociative algebra if the following identities hold on it:
\begin{gather}
(x_1\succ x_2+x_1\prec x_2)\succ x_3 = x_1\succ (x_2 \succ x_3), \label{ident:preAs1} \\
(x_1\succ x_2)\prec x_3=x_1\succ(x_2\prec x_3), \label{ident:preAs2} \\
x_1\prec(x_2\succ x_3+x_2\prec x_3)=(x_1\prec x_2)\prec x_3. \label{ident:preAs3}
\end{gather}

In \cite{Dialg99}, free preassociative algebra in the terms of forests was constructed.

Let $A$ be a preassociative algebra, $a_1,a_2,\ldots,a_k\in A$. Introduce
\begin{gather*}
(a_1,a_2,\ldots,a_k)_{\succ} = a_1\succ (a_2\succ(\ldots \succ (a_{k-2}\succ (a_{k-1}\succ a_k))\ldots)),\\
(a_1,a_2,\ldots,a_k)_{\prec} = ((\ldots((a_1\prec a_2)\prec a_3)\prec\ldots )\prec a_{k-1})\prec a_k.
\end{gather*}

{\bf Statement}.
The following relation holds in any preassociative algebra $A$:
\begin{equation}\label{eq:big-preas}
((b_1,b_2,\ldots,b_k,a)_{\succ},c_1,c_2,\ldots,c_l)_\prec
 = (b_1,b_2,\ldots,b_k,(a,c_1,c_2,\ldots,c_l)_{\prec})_{\succ},
\end{equation}
for any natural numbers $k,l$ and
$a$, $b_1,\ldots,b_k$, $c_1,\ldots,c_l\in A$.

{\sc Proof}.
We will prove the statement by induction on $k$. Let $k = 1$.
Now we will proceed the induction on $l$.
For $l=1$ the statement follows from \eqref{ident:preAs2}.
Assume that we have proved the formula for all numbers less $l$.
The inductive step on $l$ follows from the eqialities
\begin{multline*}
((b_1,a)_\succ,c_1,c_2,\ldots,c_l)_\prec
 = (\ldots (((b_1\succ a)\prec c_1)\prec c_2)\prec \ldots )\prec c_l \\
 = (\ldots((b_1\succ (a\prec c_1))\prec c_2)\prec \ldots )\prec c_l
 = ((b_1,a\prec c_1)_\succ,c_2,c_3,\ldots,c_l)_\prec \\
 = (b_1, (a\prec c_1,c_2,c_3,\ldots,c_l)_\prec )_\succ
 = (b_1, (a,c_1,c_2,\ldots,c_l)_{\prec})_{\succ}.
\end{multline*}

Assuming the statement is true for all natural numbers less $k$,
the following equalities prove the inductive step:
\begin{multline*}
((b_1,b_2,\ldots,b_k,a)_{\succ},c_1,c_2,\ldots,c_l)  \\
 = ((b_1,b_2,\ldots,b_{k-1},b_k\succ a)_{\succ},c_1,c_2,\ldots,c_l) \\
 = (b_1,b_2,\ldots,b_{k-1}, (b_k\succ a,c_1,c_2,\ldots,c_l)_\prec)_\succ \\
 = (b_1,b_2,\ldots,b_{k-1}, ((b_k,a)_\succ,c_1,c_2,\ldots,c_l)_\prec)_\succ \\
 = (b_1,b_2,\ldots,b_{k-1}, (b_k, (a,c_1,c_2,\ldots,c_l)_\prec)_\succ)_\succ \\
 = (b_1,b_2,\ldots,b_k,(a,c_1,c_2,\ldots,c_l)_{\prec})_{\succ}.
\end{multline*}

\subsection{Postassociative algebra}

Postassociative algebra is a linear space with three bilinear operations
$\cdot$, $\succ$, $\prec$ satisfying 7 identities:
\begin{equation}\label{id:PostAs}
\begin{gathered}
(x \prec y) \prec z = x \prec (y \succ z + y \prec z + y\cdot z), \quad
(x \succ y) \prec z = x \succ (y \prec z), \\
(x \succ y + y\succ x + x\cdot y) \succ z = x \succ (y \succ z), \\
x\succ (y\cdot z) = (x\succ y) \cdot z,\quad
(x \prec y) \cdot z = x \cdot (y \succ z), \\
(x \cdot y) \prec z = x \cdot (y \prec z), \quad
(x \cdot y) \cdot z = x \cdot (y \cdot z).
\end{gathered}
\end{equation}

{\bf Example 4} \cite{Trialg01}.
Let $K = \Bbbk[y_1,y_2,\ldots,y_k,\ldots]$ be polynomial algebra on countable number of variables $y_k$.
By induction, define the operations $\cdot$, $\succ$, $\prec$ on the augmentation ideal $I\triangleleft K$:
$$
y_{k}\omega \cdot y_{k'}\omega'
 = y_{k+k'}(\omega*\omega'), \quad
y_{k}\omega \succ y_{k'}\omega'
 = y_{k'}(y_k\omega*\omega'), \quad
y_{k}\omega \prec y_{k'}\omega'
 = y_{k}(\omega*y_{k'}\omega'),
$$
where $a*b = a\succ b+ a\prec b+a\cdot b$.
The space $I$ with the defined operations is postassociative algebra.

In \cite{Trialg01}, free postassociative algebra in the terms of trees was constructed.

\subsection{Embedding of Loday algebras in RB-algebras}

The common definition of the varieties of pre- and post-$\Var$-algebra
for a variety $\Var$ could be found in \cite{BaiGuoNi10} or \cite{GubKol2014}.

Given an associative algebra $B$ with RB-operator $R$ of zero weight,
the space $B$ with respect to the operations
\begin{equation}\label{eq:RB-PreOperations}
x\succ y = R(x)y, \quad
x\prec y = xR(y)
\end{equation}
is a preassociative algebra.

For any RB-operator $R$ on $B$ of weight $\lambda\neq0$,
we have that an operator $R' = \frac{1}{\lambda}R$
is an RB-operator of unit weight.
The space $B$ with the operations $x\cdot y = xy$ and
\eqref{eq:RB-PreOperations} defined for $R'$
is a postassociative algebra.
Denote the constructed pre- and postassociative algebras as $B^{(R)}_{\lambda}$.
For short, we will denote $B^{(R)}_{0}$ as $B^{(R)}$.

Given a preassociative algebra $\langle C,\succ,\prec\rangle$,
universal enveloping associative RB-algebra $U$ of $C$
is an universal algebra in the class of all associative RB-algebras
of zero weight such that there exists injective homomorphism from $C$ to $U^{(R)}$.
Analogously universal enveloping associative RB-algebra
of a postassociative algebra is defined. The common denotation
of universal enveloping of pre- or postassociative algebra: $U_{\RB}(C)$.

{\bf Theorem 1} \cite{GubKol2013}.
a) Any pre-$\Var$-algebra could be embedded into its universal enveloping
RB-algebra of the variety $\Var$ and zero weight.

b) Any post-$\Var$-algebra could be embedded into its universal enveloping
RB-algebra of the variety $\Var$ and nonzero weight.

Based on Theorem 1, we have the natural question:
What does a linear base of universal enveloping RB-algebra
of a pre- or post-$\Var$-algebra look like for an arbitrary variety $\Var$?
In the case of associative pre- and postalgebras, the question
appeared in \cite{Guo2011}.
The current article is devoted to answer the question
in the associative case.

In the article, the common method to construct universal enveloping is the following.
Let $X$ be a linear base of a preassociative algebra $K$.
We find a base of universal enveloping $U_{\RB}(K)$
as the special subset $E$ of the standard base of $\RB\As\langle X\rangle$
closed under the action of RB-operator.
By induction, we define a product $*$
on the linear span of $E$ and prove its associativity.
Finally, we state universality of the algebra $\Bbbk E$.

In the case of postassociative algebras,
as it was mentioned above, we will consider
universal enveloping associative RB-algebra of unit weight.

\section{Universal enveloping Rota---Baxter algebra of pre\-associative algebra}

In the paragraph, we will construct universal enveloping RB-algebra
of arbitrary preassociative algebra $\langle C,\succ,\prec\rangle$.
Let $B$ be a linear base of $C$.

{\sc Definition}.
An element $v = R(x)$ of the standard base $\RB\As\langle B\rangle$ is

\noindent--- left-good if $v$ is not of the form $R(b)$ or $R(bR(y))$ for $b\in B$,

\noindent--- right-good $v$ is not of the form $R(b)$ or $R(R(y)b)$ for $b\in B$,

\noindent--- good if $v$ is left- and right-good simultaneously,

\noindent--- semigood if $v$ is good or of the form
$$
v = R(a_1 R( R(a_3 R(R(\ldots R(a_{2k-1} R(R(y)a_{2k}) )  \ldots )a_4) )a_2 ) ),
$$
where $R(y)$ is good, $a_1,a_2,\ldots,a_{2k}\in B$.

Remark that any semigood element is right-good and
any element from the standard base of the form
$v = R(x)$, $x\not\in B$, is left-good or right-good.

Let us consider an associative algebra
$A = \As\langle B|(a\prec b)c - a(b\succ c),\,a,b,c\in B\rangle$
with the multiplication $\cdot$.
Here the expressions $a\prec b$ and $a\succ b$, $a,b\in B$,
equal the results of the products in the preassociative algebra $C$.
As $B$ is the linear base of $C$,
the expressions are linear combinations of the elements of $B$.

Due to Gr\"{o}bner---Shirshov theory, namely
diamond lemma for associative algebras \cite{Bergman,Bokut},
there exists such set $E_0\subset S(B)$ that $\bar{E}_0$ ---
the image of $E_0$ under the factorization of $\As\langle B\rangle$ by the ideal
$\langle (a\prec b)c - a(b\succ c),\,a,b,c\in B\rangle$
--- is a base of $A$; moreover,
for any decomposition of an element $v\in E_0$ into a concatenation
$v = v_1 v_2$ we have $v_1,v_2\in E_0$.

Let us define by induction $Envelope$-words (shortly $E$-words),
a subset of the standard base of $\RB\As\langle B\rangle$:

1) elements of $E_0$ are $E$-words of the type 1;

2) given $E$-word $u$, we define $R(u)$ as an $E$-word of the type 2;

3) the word
$$
v = u_0 R(v_1) u_1 R(v_2) u_2 \ldots u_{k-1} R(v_k) u_k,
\quad \deg v\geq 2,\,\deg_R(v)\geq1,
$$
is an $E$-word of the type 3, if
$u_1,\ldots, u_{k-1}\in E_0$, $u_0,u_k\in E_0\cup\{\emptyset\}$,
$v_1,\ldots,v_{k}$ are $E$-words,
$R(v_2),\ldots,R(v_{k-1})$ are semigood elements of the standard base,
$R(v_1)$ is right-good.
Given $u_0\neq\emptyset$, the element $R(v_1)$ is semigood.
Given $u_k\neq\emptyset$, the element $R(v_k)$ is semigood,
else $R(v_k)$ is left-good.

{\bf Theorem 2}.
The set of all $E$-words forms a linear base of
universal enveloping associative RB-algebra of $C$.

{\bf Lemma 1}.
Let $D$ denote a linear span of all $E$-words.
One can define such bilinear operation $*$ on the space $D$ that
($k$--$l$ denotes below the condition on the product $v*u$,
where $v$ is an $E$-word of the type $k$ and $u$ is an $E$-word of the type $l$.)

1--1: given $v,u\in E_0$, we have
\begin{equation}\label{def:1--1}
v*u = v\cdot u.
\end{equation}

1--2: given $v = w'a\in E_0$, $a\in B$, $w'\in E_0\cup \{\emptyset\}$, an $E$-word $u = R(p)$ of the type 2, we have
\begin{equation}\label{def:1--2}
v*u = \begin{cases}
   w'\cdot (a\prec b), & p = b \in B, \\
  (w'\cdot(a\prec b))R(x)-wR(R(b)*x), & p = bR(x), b\in B, \\
   wR(p), & R(p)\ \mbox{is left-good}.
               \end{cases}
\end{equation}

1--3: given $v=w=w'a\in E_0$, $a\in B$, $w'\in E_0\cup \{\emptyset\}$, an $E$-word $u$ of the type 3, we have
\begin{equation}\label{def:1--3}
v*u
= \begin{cases}
  (w\cdot u_0)R(x)u', & u = u_0 R(x)u',\,u_0\in E_0, \\
   (w'\cdot (a\prec b))*(R(y)*u') & u = R(bR(y))u',\,b\in B,\\
   \hfill - wR(R(b)*y)*u' & R(bR(y))\ \mbox{is not semigood}, \\
   wR(x)u', & u = R(x)u',\,R(x)\ \mbox{is semigood}.
      \end{cases}
\end{equation}

3--1: given $u=w=aw'\in E_0$, $a\in B$, $w'\in E_0\cup \{\emptyset\}$, an $E$-word $v$ of the type 3, we have
\begin{equation}\label{def:3--1}
v*u
= \begin{cases}
  v'R(x)(v_0\cdot w), & v = v'R(x)v_0,\,v_0\in E_0, \\
  (v'*R(y))*((b\succ a)\cdot w') & v = v'R(R(y)b),\,b\in B, \\
 \hfill  - v'*R(y*R(b))w, & \\
  v'R(x)w, & v = v'R(x),\,R(x)\ \mbox{is good}.
      \end{cases}
\end{equation}

2--2: given $E$-words $v = R(v')$, $u = R(u')$ of the type 2, we have
\begin{equation}\label{def:2--2}
v*u = R(v') * R(u') = R( R(v')*u' + v'*R(u')).
\end{equation}

2--3: given an $E$-word $v = R(x)$ of the type 2, an $E$-word $u$ of the type 3, we have
\begin{equation}\label{def:2--3}
v*u
= \begin{cases}
  ((a\succ b)\cdot u_0)R(y)u', & x=a,u = bu_0 R(y)u',\,u_0\in E_0,\,a,b\in B, \\
    R(z)((a\succ b)\cdot u_0)R(y)u' & x=R(z)a,\,u = bu_0 R(y)u',\,u_0\in E_0,a,b\in B, \\
    \hfill -R(z*R(a))u, &  \\
    R(x)u, & u = bu',\,b\in B,\,R(x)\ \mbox{is right-good}, \\
   R( R(x)*y + x*R(y) )u', & u = R(y)u'.
      \end{cases}
\end{equation}

3--3: given $E$-words $v,u$ of the type 3, we have
\begin{equation}\label{def:3--3}
v*u = \begin{cases}
  v' R(R(x)*y + x*R(y) )u', & v = v' R(x),\,u = R(y)u', \\
  v'R(x)(w*u), & v = v'R(x)w,\,u = R(y)u',w\in E_0, \\
  (v*w)R(y)u', & v = v'R(x),\,u = wR(y)u',w\in E_0, \\
  v'R(x)(v_0\cdot u_0)R(y)u', & v = v'R(x)v_0,\,u = u_0R(y)u',\,v_0,u_0\in E_0.
      \end{cases}
\end{equation}

The following conditions are also satisfied:

L1) Given left-good $E$-word $R(x)$, $b\in B$,
and $R(b)*x = \sum u_i$, where $u_i$ is an $E$-word,
we have that $R(u_i)$ is left-good for every $i$.

L2) Given right-good $E$-word $R(y)$,
any $E$-word $x$, and $R(x)*y = \sum u_i$,
where $u_i$ is an $E$-word,
we have that $R(u_i)$ is right-good for every $i$.

L3) Given $E$-word $x\not\in B$, $E$-word $u$ of the type 1 or 3,
and $R(x)*u = \sum u_i$, where $u_i$ is an $E$-word,
we have that $u_i$ is an $E$-word of the type 3
and begins with an $R$-letter, i.e., has a view $u_i = R(x_i)u_i'$, for every $i$.

L4) Given an $E$-word $u = u'a$, $a\in B$, of the type 1 or 3,
any $E$-word $v$, and $v*u = \sum u_i$, where $u_i$ is an $E$-word,
we have $u_i = u_i'a_i$, $a_i\in B$.

L5) Given $E$-word $u = aw$, $a\in B$, $w\in E_0\cup\{\emptyset\}$, of the type 1,
any $E$-word $v$, and $v*u = \sum u_i$, where $u_i$ is an $E$-word,
we have that $u_i$ is one of the following forms:
\begin{gather*}
w'\cdot aw, \quad
(c\succ a)\cdot w,\quad
v_i'R(x_i)((c_1,c_2,\ldots,c_k,a)_\succ\cdot w),
\end{gather*}
for $c,c_1,c_2,\ldots,c_k\in B$, $k\geq0$, $w'\in E_0$.
Moreover, the number, view of summands and values of $v_i'R(x_i)$ and $c_l$
depend only on $v$.

L6) Given $E$-word $u = awR(p)u'$, $a\in B$, $w\in E_0\cup\{\emptyset\}$, of the type 3,
any $E$-word $v$, and $v*u = \sum u_i$, where $u_i$ is an $E$-word,
we have that $u_i$ is one of the following forms:
\begin{gather*}
(w'\cdot aw) R(p)u',\quad
((c\succ a)\cdot w)R(p)u',\quad
v_i'R(x_i)((c_1,c_2,\ldots,c_k,a)_\succ\cdot w) R(p)u',
\end{gather*}
where $c_1,c_2,\ldots,c_k\in B$, $k\geq0$, $w'\in E_0$.
Moreover, the number, view of summands and values of $v_i'R(x_i)$ and $c_l$
depend only on $v$.

L7) Given $E$-word $u = R(s)awR(t)u'$, $a\in B$, $w\in E_0\cup\{\emptyset\}$, of the type 3,
in which empty word could stay instead of $R(t)u'$, any $E$-word $v$, and
$v*u = \sum u_i$, where $u_i$ is an $E$-word of the form
$$
v_i'R(s_i)((c_1,c_2,\ldots,c_k,a)_\succ\cdot w) R(t)u',
$$
$c_1,c_2,\ldots,c_k\in B$, $k\geq0$.
Moreover, the number, view of summands and values of $v_i'R(s_i)$ and $c_l$
depend only on $v$ and $s$.

The relations for the products of $E$-words
of types 2--1, 3--2 and the conditions R1--R7
are defined analogously to the products
of $E$-words of types 1--2, 2--3 and the conditions
L1--L7 by the inversion of letters, multipliers
and operations symbols, wherein
$\succ$ and $\prec$ turn in each other.

{\sc Proof}.
Let us define the operation $*$ with the prescribed conditions
for $E$-words $v,u$ by induction on $r = \deg_R(v)+\deg_R(u)$.
For $r = 0$ define $v*u = v\cdot u$, $v,u\in E_0$,
which satisfies the condition 1--1. It is easy to
see that the conditions L5 and R5 are also fulfilled,
all others are true because of $r=0$.

For $r = 1$, let us define $v*u$ by induction on
$d = \deg(v)+\deg(u)$. For $d = 2$,
$v=a\in B$, $u=R(w)$, $w\in E_0$, we define
\begin{equation}\label{def*:r=0d=2}
v*u = a*R(w) = \begin{cases}
a\prec b, & w = b\in B, \\
aR(w), & w\in E_0\setminus B,
\end{cases}
\end{equation}
which satisfies the condition 1--2.
We define the product $R(w)*a$ analogously
to \eqref{def*:r=0d=2} up to the inversion.
The cases 1--1 and 2--2
are not realizable because of $r = 1$.
It is correct to write $aR(w)$, as
the element $R(w)$ is left-good for $w\in E_0\setminus B$.

For $r = 1$, $d > 2$, we define $v*u$ for
pairs of $E$-words of types 1--2, 1--3;
the definition for types 2--1, 3--1 is analogous
up to the inversion.
The cases 2--3, 3--2, and 3--3 do not appear for $r = 1$.

Let $v = w_1 = w'a\in E_0$, $a\in B$, $u = R(w_2)$, $w_1,w_2\in E_0$, define
$$
v*u  = \begin{cases}
w'\cdot(a\prec b), & w_2 = b\in B, \\
w_1R(w_2), & w_2\in E_0\setminus B,
         \end{cases}
$$
what is consistent with the condition 1--2. The notation $w_1R(w_2)$
is correct because the element $R(w_2)$ is left-good
for $w_2\in E_0\setminus B$.

Let $v\in E_0$, $u$ be an $E$-word of the type 3, define
$$
v*u = \begin{cases}
(v\cdot u_0)R(w)u', & u = u_0R(w)u',u_0\in E_0,u'\in E_0\cup\{\emptyset\},w\in E_0\setminus B, \\
vR(w)u', & u = R(w)u',u'\in E_0,w\in E_0\setminus B,
         \end{cases}
$$
what is consistent with the condition 1--2.
The notations are correct by the same arguments as above.

For $r=1$, clarify that the conditions L1--L7
(and analogously R1--R7) hold. Indeed,
L1) $R(x)$ is left-good, $\deg_R(x) = 0$,
so $x \in E_0$, $\deg(x)\geq2$, $x=ax'$, $a\in B$.
Hence, $R(b)*x = (b\succ a)\cdot x'$ is a linear combination
of the elements from $E_0\setminus B$.
L3) Given left-good $E$-word $R(x)$,
we have $x\in E_0\setminus B$
and $R(x)*w = R(x)w$ for $w\in E_0$.
L2) Given right-good $E$-word $R(y)$, $y\in E_0$,
we have $y = ay'\in E_0\setminus B$, $a\in B$. So
$R(x)*y$ equal $R(x)y$ for $x\in E_0\setminus B$
or $(x\succ a)\cdot y'$ for $x\in B$
is a linear combination of the elements from $E_0\setminus B$.
It is easy to check the conditions L4--L7 (R1--R7).

Suppose that the product $v*u$ is yet defined for
all pairs of $E$-words $v,u$ such that
$\deg_R(v)+\deg_R(u)<r$, $r\geq2$, and all conditions
of the statement on $*$ are satisfied.
Let us define $*$ on $E$-words $v,u$ with
$\deg_R(v)+\deg_R(u) = r$
by induction on $d = \deg(v)+\deg(u)$.

For $d = 2$, consider the cases 1--2 and 2--2
(the case 2--1 is analogous to 1--2).

1--2: let $v = a\in B$, $u = R(p)$ be an $E$-word of the type 2,
$\deg_R(u)=r\geq 2$, define $v*u$ by \eqref{def:1--2}.

2--2: let $v = R(v')$, $u = R(u')$ be an $E$-word of the type 2,
define $v*u$ by \eqref{def:2--2}.

The products $R(b)*x$, $R(v')*u'$, $v'*R(u')$ in
\eqref{def:1--2}, \eqref{def:2--2}
are defined by induction on $r$.
The element $p = bR(x)$ is an $E$-word of the type 3,
so $R(x)$ is left-good and, therefore,
the concatenation $(a\prec b)R(x)$ is correct.
By the condition L1 for $R(b)*x$
holding by inductive assumption,
the notation $aR(R(b)*x)$ is correct.
The conditions L2--L5 (and also R2--R5)
in the cases $d=2$, 2--1 and 2--2 hold;
the conditions L1 and R1
are realizable only in the case 2--2,
hence, they are also fulfilled;
the conditions L6, L7 (as R6, R7)
are totally not realizable, so they hold.

For $d > 2$, define the product $v*u$ for
$E$-words pairs of the following cases: 1--2, 2--3, 3--3, 1--3, 3--1;
the products for the cases 2--1 and 3--2 are defined
analogously up to the inversion.

1--2: let $v = w = w'a\in E_0$, $a\in B$,
$u = R(p)$ be an $E$-word of the type 2,
$\deg_R(u)=r\geq 2$, define $v*u$ by \eqref{def:1--2}.

The product $R(b)*x$ in \eqref{def:1--2}
is defined by induction on $r$.
The element $p = bR(x)$ is $E$-word of the type 3,
so $R(x)$ is left-good and the concatenation
$(w'\cdot(a\prec b))R(x)$ is correct.
By the condition L1 holding for $R(b)*x$
by inductive assumption,
the notation $wR(R(b)*x)$ is correct.
The definition of $v*u$ is consistent with L1--L7 (R1--R7).

2--3: let $v = R(x)$ be an $E$-word of the type 2,
$u$ be an $E$-word of the type 3, define $v*u$ by \eqref{def:2--3}.

The notations $R(z*R(a))u$, $R( R(x)*y + x*R(y) )u'$
in \eqref{def:2--3} are correct by the conditions R1, L2, and R3;
the products are defined by induction on $r$.
The definition of $v*u$ is consistent with L1--L7 (R1--R7).

3--3: let $v,u$ be $E$-words of the type 3,
define $v*u$ by \eqref{def:3--3}.

The correctness of the products in \eqref{def:3--3}
follows in the first case by the conditions L2, L3, R2, R3,
in the second and third --- by L4, R4.
The definition of $v*u$ is consistent with L1--L7 (R1--R7).

1--3: let $v=w=w'a\in E_0$, $a\in B$,
$u$ be an $E$-word of the type 3,
define $v*u$ by \eqref{def:1--3}.

3--1: let $u=w=aw'\in E_0$, $a\in B$, $v$ be an $E$-word of the type 3,
define $v*u$ by \eqref{def:3--1}.

The conditions L1 and R1 provide that the notations $wR(R(b)*y)$,
$R(y*R(b))w$ in \eqref{def:1--3}, \eqref{def:3--1} are correct.
The definition of the product in the case 1--3 for $u = R(bR(y))u'$,
where $R(bR(y))$ is not semigood, is reduced to the cases 3--1 and 3--3.
In the last variant of the case 3--1,
the product is expressed by the one from the case 1--3.
We have to prove that the process of computation $v*u$ in the cases
1--3 and 3--1 is finite.
Suppose we have $v = w=w'b\in E_0$, $b\in B$,
$u = R(bR(y))u'$, and $R(bR(y))$ is not semigood, i.e., has the form
$$
u = R(a_1 R( R(a_3 R(R(\ldots R(a_{2k+1} R(y))  \ldots )a_4) )a_2 ) )cw'',
$$
where $R(y)$ is good, $a_1,a_2,\ldots,a_{2k+1},c\in B$.
Provided $t-p$ is a sum of $E$-words or products
with summary $R$-degree less than $r$,
we will denote it by $t\equiv p$. Write
\begin{multline}\label{eq:corr1}
\allowdisplaybreaks
v*u \equiv
 - vR(R(a_1)*(R(a_3 R(R(\ldots R(a_{2k+1} R(y))  \ldots )a_4) )a_2) )*cw'' \\
 \equiv - vR( R( a_1*R(a_3 R(R(\ldots R(a_{2k+1} R(y))  \ldots )a_4) ))a_2 )*cw'' \\
 \equiv   vR( R( a_1 R( R(a_3)*(R( a_5 R(R(\ldots R(a_{2k+1} R(y))  \ldots )a_6 )))a_4 ) ) )a_2 )*cw'' \\
 \equiv - vR( R( a_1 R( R(a_3*R(a_5 R(R(\ldots R(a_{2k+1} R(y))  \ldots )a_6 ))) a_4 ) )a_2 )*cw''.
\end{multline}
Continuing on and rewriting analogously the product into the action of the central $R$-letter,
finally we will have the expression
\begin{equation}\label{eq:corr2}
v*u \equiv
  - vR( R( a_1 R( R(a_3 R(\ldots R( R( a_{2k-1}R( R(a_{2k+1})*y) )a_{2k})\ldots ) a_4 ) )a_2 )*cw'',
\end{equation}
in which $R(R(a_{2k+1})*y)$ is left-good by the condition L1 and right-good by L2.
Let
\begin{equation}\label{eq:corr3}
R(a_{2k+1})*y = \sum z_i,
\end{equation}
for good $E$-words $R(z_i)$.
We will prove that the definition of $v*u$ is correct and has no cycles for any $i$.
By the definition of the product for types 3--1, we have
\begin{multline}\label{eq:corr4}
vR( R( a_1 R( R(a_3 R(\ldots R( R( a_{2k-1}R( z_i ) )a_{2k})\ldots )) a_4 ) )a_2 )*cw'' \\
\equiv
- v*R( a_1 R( R(a_3 R(\ldots R( R( a_{2k-1}R( z_i ) )a_{2k})\ldots )) a_4 )*R(a_2))cw'' \\
\equiv
- v*R( a_1 R( R( R(a_3 R(\ldots R( R( a_{2k-1}R( z_i ) )a_{2k})\ldots )) a_4 )*a_2) )cw'' \\
\equiv
  v*R( a_1 R( R( R(a_3 R(\ldots R( R( a_{2k-1}R( z_i ) )a_{2k})\ldots )) a_4 )*a_2) )cw'' \\
\equiv
- v*R( a_1 R( R( (a_3 R(\ldots R( R( a_{2k-1}R( z_i ) )a_{2k})\ldots ))*R(a_4) )a_2) )cw'' \\
\equiv \ldots \equiv
- v*R( a_1 R( R( (a_3 R(\ldots R( a_{2k-1}R( z_i )*R(a_{2k})\ldots )) a_4 ) )a_2) )cw'' \\
= - vR( a_1 R( R( (a_3 R(\ldots R( a_{2k-1}R( R(z_i)a_{2k}+z_i*R(a_{2k}) )\ldots )) a_4 ) )a_2) )cw'',
\end{multline}
$R(z_i*R(a_{2k}))$ is good by R1 and R2, so the last equality is true.

The product $v*u$ of $E$-words of types 1--3 and 3--1
satisfies the conditions L1--L7 (R1--R7).

{\bf Lemma 2}.
The space $D$ with the operations $*$, $R$ is an RB-algebra.

{\sc Proof}.
It follows from \eqref{def:2--2}.

{\bf Lemma 3}.
The relations
$R(a)*b = a\succ b$, $a*R(b) = a\prec b$ hold in $D$ for every $a,b\in B$.

{\sc Proof}.
It follows from Lemma 1, the first case of \eqref{def:1--2}, and
analogous relation of 2--1.

{\bf Lemma 4}.
Given any $E$-word $w\in E_0$, $a,b\in B$,
the equality $(w,a,R(b)) =  (R(b),a,w) = 0$
is true on $D$.

{\sc Proof}.
Let us define a map
$\dashv\colon \As\langle B\rangle\otimes \Bbbk B\to \As\langle B\rangle$
on the base as follows:
\begin{equation}\label{eq:dashv1}
wa\dashv b = w(a\prec b),
\end{equation}
where $w\in S(B)\cup\emptyset$, $a,b\in B$,
and $a\prec b$ is a product in $C$
and is equal to a linear combination of elements of $B$.
Let us check that the map $\dashv$
could be induced on the algebra $A = \As\langle B|(a\prec b)c - a(b\succ c),\,a,b,c\in B\rangle$.
Let
$$
I = \langle (a\prec b)c - a(b\succ c),\,a,b,c\in B\rangle\lhd \As\langle B\rangle.
$$
From \eqref{ident:preAs2} and linearity,
we have $(b\succ c)\dashv d = b\succ (c\prec d)$, $b,c,d\in B$.
Applying this equality, we obtain
\begin{multline*}
(wa(b\succ c) - w(a\prec b)c)\dashv d
 = wa((b\succ c)\dashv d) - w(a\prec b)(c\prec d) \\
 = wa(b\succ (c\prec d)) - w(a\prec b)(c\prec d)\in I,
\end{multline*}
$$
(a(b\succ c)we - (a\prec b)cwe)\dashv d
 = a(b\succ c)w(e\prec d) - (a\prec b)cw(e\prec d) \in I
$$
for $a,b,c,d\in B,\ w\in E_0\cup\{\emptyset\}$.

Hence, the induced map
$\dashv\colon A\otimes \Bbbk B\to A$ is well-defined.
Notice that
\begin{equation}\label{eq:dashv2}
(wa)*R(b) = w\cdot (a\prec b) = wa\dashv b + I,\quad
a,b\in B,\ w\in E_0\cup\{\emptyset\},
\end{equation}
where $wa\in E_0$, $a\in B$, $w\in E_0\cup\{\emptyset\}$.
By \eqref{eq:dashv1} and \eqref{eq:dashv2}, conclude
$$
(w,a,R(b))
 = (w\cdot a)\dashv b - w\cdot (a\prec b) = 0.
$$
The proof of the equality $(R(b),a,w) = 0$ is analogous.

{\bf Lemma 5}.
The operation $*$ on $D$ is associative.

{\sc Proof}. Given $E$-words $x,y,z$, let us prove associativity
$$
(x,y,z) = (x*y)*z - x*(y*z) = 0
$$
by inductions on two parameters:
at first, on summary $R$-degree $r$ of the triple $x,y,z$,
at second, on summary degree $d$ of $x,y,z$.

For $r = 0$, we have $(x,y,z) = (x\cdot y)\cdot z - x\cdot (y\cdot z) = 0$, $x,y,z\in E_0$,
as the product $\cdot$ is associative in the algebra $A$.

Let $r>0$ and suppose that associativity
for all triples of $E$-words with the less summary $R$-degree is proven.

We prove the statement for the triples $x,y,z$,
in which $y$ is an $E$-word of the type 1 or 3, $d$ is any.
Let $y$ be an $E$-word of the type 1. Consider the case $y = a$, $a\in B$.
By the condition L5, the product $x*y$ is a sum $\sum\limits_{i\in I} u_i$,
where $u_i$ has one of the following forms:
$$
w'\cdot a,\quad
e\succ a,\quad
v_i'R(x_i)(c_1,c_2,\ldots,c_k,a)_\succ,
$$
$e,c_1,\ldots,c_k\in B$.

By R5, the product $u_i*z$ equals to $\sum\limits_{j\in J}t_{ij}$,
where $t_{ij}$ is one of following forms:
(The antepenultimate case is written by Lemma 4):
\begin{equation}\label{eq:middle-letter-l}
\begin{gathered}
(w'\cdot a)\cdot w'', \quad
(e\succ a)\cdot w'',\quad
v_i'R(x_i)((c_1,c_2,\ldots,c_k,a)_\succ\cdot w''),\\
w'\cdot (a\prec f),\quad
(e\succ a)\prec f,\quad
v_i'R(x_i)((c_1,c_2,\ldots,c_k,a)_\succ\prec f),\\
(w'\cdot(a,d_1,d_2,\ldots,d_l)_\prec )R(y_i)u_i',\quad
(e\succ a,d_1,d_2,\ldots,d_l)_\prec R(y_i)u_i',\\
v_i'R(x_i)((c_1,c_2,\ldots,c_k,a)_\succ,d_1,d_2,\ldots,d_l)_\prec R(y_i)u_i',
\end{gathered}
\end{equation}
$f,d_1,\ldots,d_l\in B$.

By the conditions R5 and L5, the product $x*(y*z)$
is a sum with the same indexing sets $I,J$ and summands of the forms
\begin{equation}\label{eq:middle-letter-r}
\begin{gathered}
w'\cdot (a\cdot w''), \quad
(e\succ a)\cdot w'',\quad
v_i'R(x_i)((c_1,c_2,\ldots,c_k,a)_\succ\cdot w''),\\
w'\cdot (a\prec f),\quad
e\succ (a \prec f),\quad
v_i'R(x_i)(c_1,c_2,\ldots,c_k,a\prec f)_\succ,\\
(w'\cdot(a,d_1,d_2,\ldots,d_l)_\prec )R(y_i)u_i',\quad
(e\succ (a,d_1,d_2,\ldots,d_l)_\prec)R(y_i)u_i',\\
v_i'R(x_i)(c_1,c_2,\ldots,c_k,(a,d_1,d_2,\ldots,d_l)_\prec)_\succ R(y_i)u_i'.
\end{gathered}
\end{equation}

The corresponding summands in \eqref{eq:middle-letter-l} and \eqref{eq:middle-letter-r}
either coincide or are equal by associativity in the algebra $A$,
the equality \eqref{ident:preAs2}, and Statement.

For $y = aw$, $a\in B$, $w\in E_0$, the proof
of associativity is analogous to the case $y = a$.

Let us consider the case when $y$ is an $E$-word of the type 3.
Given $y = aR(p)u'$, $a\in B$, $u'\neq\emptyset$,
by \eqref{def:3--1}, \eqref{def:2--3}, \eqref{def:3--3}, and
the conditions L6 and R4, we have
$(x*y)*z = \sum u_i$, where $u_i$ is one of the forms:
\begin{equation}\label{eq:middle-3-y-type1}
(w'\cdot a) R(p)(u'*z),\quad
(c\succ a)R(p)(u'*z),\quad
u_i'R(x_i)(c_1,c_2,\ldots,c_k,a)_\succ R(p)(u'*z).
\end{equation}
By \eqref{def:1--3}, \eqref{def:2--3}, \eqref{def:3--3},
and the conditions L6, R4,
$x*(y*z) = \sum u_i$, where $u_i$ is one of the forms
listed in \eqref{eq:middle-3-y-type1}.
It proves associativity for the triple $x,y,z$.

Let $y = aR(p)$, $a\in B$. By the conditions L6, R7, R3,
$(x*y)*z = \sum u_i$, where $u_i$ has one of the following forms:
\begin{equation}\label{eq:middle-3-y-type2l}
\begin{gathered}
((w'\cdot a),d_1,\ldots,d_l)_\prec R(p_i), \quad
((e\succ a),d_1,\ldots,d_l)_\prec R(p_i),\\
v_i'R(x_i)((c_1,c_2,\ldots,c_k,a)_\succ,d_1,\ldots,d_l)_\prec R(p_i),
\end{gathered}
\end{equation}
$e,c_1,\ldots,c_k,d_1,\ldots,d_l\in B$, $k,l\geq0$.

By the conditions L6, R7, R4,
$x*(y*z) = \sum u_i$, where $u_i$ has one of the following forms:
\begin{equation}\label{eq:middle-3-y-type2r}
\begin{gathered}
(w'\cdot (a,d_1,\ldots,d_l)_\prec)R(p_i), \quad
(e\succ (a,d_1,\ldots,d_l)_\prec)R(p_i),\\
v_i'R(x_i)(c_1,c_2,\ldots,c_k,(a,d_1,\ldots,d_l)_\prec)_\succ R(p_i),
\end{gathered}
\end{equation}
The corresponding summands in \eqref{eq:middle-3-y-type2l} and \eqref{eq:middle-3-y-type2r}
are equal by Statement and Lemma 4.

Let $y = R(s)aR(t)$, $a\in B$. By the conditions L7, R7, R3,
$(x*y)*z = \sum u_i$, where $u_i$ has one of the following forms:
$$
v_i'R(s_i)((c_1,c_2,\ldots,c_k,a)_\succ,d_1,\ldots,d_l)_\prec R(t_i)z'.
$$
By the conditions L7, R7, R4,
$x*(y*z) = \sum u_i$, where $u_i$ is one of the following forms:
$$
v_i'R(s_i)(c_1,c_2,\ldots,c_k,(a,d_1,\ldots,d_l)_\prec)_\succ R(t_i)z'.
$$
By Statement, we have associativity.

The cases $y = awR(p)u'$, $y = R(s)awR(t)u'$, where $a\in B$, $w\in E_0$,
are analogous to ones considered above.

Hence, we have to prove associativity only for triples
$x,y,z$, in which $y$ is an $E$-word of the type 2.
The definition of $*$ by Lemma 1 is symmetric
with respect to the inversion, except the cases 1--3 and 3--1, so
associativity in the triple of types $k$--2--$l$
leads to associativity in the triple $l$--2--$k$.

Prove associativity for $d = 3$ and $E$-word $y$ of the type~2.

1--2--1. Consider possible cases of $y$.
a) For $y = R(b)$, $b\in B$, we have
$$
(a*R(b))*c - a*(R(b)*c) = (a\prec b)\cdot c - a\cdot (b\succ c) = 0.
$$

b) If $y = R(u)$ is good, then
$$
(x*y)*z - x*(y*z) = aR(u)c - aR(u)c = 0.
$$

c) For $y = R(R(p)b)$, $b\in B$, we have
\begin{multline*}
(x*y)*z - x*(y*z)
 = (aR(R(p)b))*c - a*(R(R(p)b)*c) \\
 = (a*R(p))*(b\succ c) - a*R(p*R(b))c
  - a*R(p)(b\succ c) + a*R(p*R(b))c \\
  = (a,R(p),b\succ c) = 0,
\end{multline*}
the last is true by the induction on $r$.

d) For $y = R(bR(p))$, $b\in B$, we have
\begin{multline*}
(x*y)*z - x*(y*z)
 = (a*R(bR(p)))*c - a*(R(bR(p))c) \\
 = (a\prec b)R(p)*c - aR(R(b)*p)*c
  - (a\prec b)R(p)*c + aR(R(b)*p)*c = 0.
\end{multline*}

1--2--2.
a) Given left-good $y$,
associativity follows from \eqref{def:1--2}, \eqref{def:2--3}.

b) Let $y = R(b)$, $b\in B$. If $z = R(c)$, $c\in B$,
then by \eqref{ident:preAs3}
$$
(a*R(b))*R(c) - a*(R(b)*R(c)) = (a\prec b)\prec c
- a\prec(b\succ c + b\prec c) = 0.
$$
If $z = R(u)$ is left-good, then
\begin{multline*}
a*(R(b)*R(u)) = a*R(R(b)*u + bR(u)) =
a*R(R(b)*u) - a*R(R(b)*u) + (a\prec b)R(u)  \\
= (a*R(b))*R(u).
\end{multline*}
Finally, if $u = cR(p)$, $c\in B$, where $R(p)$ is left-good,
then from one hand,
\begin{multline}\label{AsL:1-R(b)-R(cR(p))}
(a*R(b))*R(cR(p)) = (a\prec b)*R(cR(p)) \\
= - (a\prec b)R(R(c)*p) + ((a\prec b)\prec c)R(p).
\end{multline}
From another hand,
\begin{multline}\label{AsR:1-R(b)-R(cR(p))}
\allowdisplaybreaks
a*(R(b)*R(cR(p)))
= a*R(R(b)*cR(p) + b*R(cRp)) \\
= a*R((b\succ c)R(p)) + a*R(-bR(R(c)*p)+(b\prec c)R(p)) \\
= -aR(R((b\succ c))*p) + (a\prec(b\succ c))R(p) \\
+ a*R( R(b)*(R(c)*p) ) - (a\prec b)R(R(c)*p) \\
- aR(R(b\prec c)*p) + (a\prec(b\prec c))R(p).
\end{multline}
Subtracting \eqref{AsR:1-R(b)-R(cR(p))} from \eqref{AsL:1-R(b)-R(cR(p))},
by \eqref{ident:preAs3} we have
\begin{multline*}
aR(R((b\succ c))*p) - a*R( R(b)*(R(c)*p) )
+ aR(R(b\prec c)*p) \\
= a*R((R(b), R(c), p )) = 0,
\end{multline*}
which is true by induction on $d$.

c) The last case, $y = R(bR(u))$, $b\in B$; $z = R(v)$.
From one hand,
\begin{multline}\label{AsL:1-R(bR(x))-R(u)}
(x*y)*z
= ( a*R(bR(u)) )*R(v) =
( -aR(R(b)*u) + (a\prec b)R(u))*R(v) \\
= -a*R( R(R(b)*u)*v + (R(b)*u)*R(v) )
+ (a\prec b)*R(R(u)*v+u*R(v)).
\end{multline}
From another hand, applying the condition L2, we have
\begin{multline}\label{AsR:1-R(bR(x))-R(u)}
x*(y*z)
= a*(R(bR(u))*R(v) ) = a*R(R(bR(u))*v+bR(u)*R(v)) \\
= a*R(R(bR(u))*v) + a*R(bR(R(u)*v+u*R(v)) \\
= a*R(R(bR(u))*v) - a*R(R(b)*(R(u)*v+u*R(v)) \\
+(a\prec b)*R(R(u)*v+u*R(v)).
\end{multline}
Subtracting \eqref{AsR:1-R(bR(x))-R(u)} from
\eqref{AsL:1-R(bR(x))-R(u)} we get the expression $a*R(\Delta)$ for
\begin{multline*}
\Delta = - R(R(b)*u)*v - (R(b)*u)*R(v) - R(bR(u))*v + R(b)*(R(u)*v+u*R(v)) \\
= - (R(b),u,R(v)) - (R(b),R(u),v) = 0,
\end{multline*}
the conclusion is true by induction on $r$.

2--2--2. Given $x = R(u)$, $y = R(v)$, $z = R(p)$, by induction we have
$$
 (R(u),R(v),R(p)) = R(  (R(u),R(v),p)
 +  (R(u),y,R(p)) +  (x,R(v),R(p)) = 0.
$$

Let $d > 3$, consider triples $x,y,z$, in which $y$ is an $E$-word of the type 2.

Show that the cases with $x$ or $z$ equal to $w\in E_0$
are reduced to the cases with $w = a\in B$. Indeed, let $x = w = w'a$, $a\in B$, $w'\in E_0$.
By the condition R4 and associativity in the triples $k$--1--$l$ and $k$--3--$l$, we have
\begin{gather*}
(x*y)*z
 = ((w'a)*y)*z
 = (w'*(a*y))*z = w'*((a*y)*z), \\
x*(y*z)
 = w'a*(y*z)
 = w'*(a*(y*z)).
\end{gather*}

Let at least one $E$-word from $x,z$ be an $E$-word of the type 3, e.g., $z$.
For $z = az'$, $a\in B$, by the condition L4,
associativity in the triples $k$--1--$l$, $k$--3--$l$,
and induction on $r$, we have
\begin{multline*}
(x*y)*z - x*(y*z)
 = (x*y)*(az') - x*(y*(az'))
 = ((x*y)*a)*z'- x*((y*a)*z') \\
 = ((x*y)*a)*z'- (x*(y*a))*z'
 = (x,y,a)*z' = 0.
\end{multline*}
If $z = R(s)az'$, $a\in B$, $y = R(p)$, then we have
\begin{multline*}
\allowdisplaybreaks
(x*y)*z - x*(y*z)
 = (x*R(p))*(R(s)az') - x*(R(p)*(R(s)az')) \\
 = ((x*R(p))*R(s))*(az') - x*((R(p)*(R(s))*az') \\
 = (x*(R(p)*R(s)))*(az') - x*((R(p)*(R(s))*az') \\
 = (x*R(R(p)*s + p*R(s)))*(az') - x*(R(R(p)*s + p*R(s))*az') \\
 = (x,R(R(p)*s + p*R(s)),az') = 0
\end{multline*}
by associativity in the triple $k$--3--$l$ and induction on $d$.

We have considered all possible cases of triples $x,y,z$ of $E$-words.
Lemma 5 is proved.

{\sc Proof of Theorem 2}.
Let us prove that the algebra $D$ is
exactly universal enveloping algebra for the preassociative algebra $C$,
i.e., is isomorphic to the algebra
$$
U_{\RB\As}(C) = \RB\As\langle B|a\succ b = R(a)b,\,a\prec b = aR(b),\,a,b\in B\rangle.
$$
By the construction, the algebra $D$ is generated by $B$.
Therefore, $D$ is a homomorphic image of a homomorphism $\varphi$ from $U_{\RB\As}(C)$.
We will prove that all basic elements of $U_{\RB\As}(C)$ are linearly expressed by $D$,
this leads to nullity of kernel of $\varphi$ and $D\cong U_{\RB\As}(C)$.

As the equality
$$
(a\prec b)c = aR(b)c = a(b\succ c)
$$
is satisfied on $U_{\RB\As}(C)$,
the space $U_{\RB\As}(C)$ is a subspace of
$\RB\As\langle B|(a\prec b)c = a(b\succ c),\,a,b,c\in B\rangle$.
It is well-known \cite{Guo2011} that algebras
$\RB\As\langle B\rangle$ and $\RB\As(\As\langle B\rangle)$ coincide.
Consider the map
$\varphi\colon B\to \RB\As(\As\langle B\rangle)/I$,
where $I$ is an RB-ideal (i.e., the ideal closed under RB-operator) generated by the set
$\{(a\prec b)c = a(b\succ c),\,a,b,c\in B\}$,
and the map $\psi$, the composition of trivial maps
$\psi_1\colon B\to A$ and $\psi_2\colon A\to\RB\As(A)$,
where as above,
$$
A = \As\langle B|(a\prec b)c - a(b\succ c),\,a,b,c\in B\rangle,\quad
\RB\As(A) = \RB\As\langle E_0|w_1 w_2 = w_1\cdot w_2\rangle.
$$
The base of $\RB\As(A)$ \cite{FardGuo08}
could be constructed by induction, it contains the elements
$$
w_1 R(u_1)w_2 R(u_2)\ldots w_k R(u_k)w_{k+1},
$$
where $w_i\in E_0$, $i=2,\ldots,k$, $w_1,w_{k+1}\in E_0\cup\{\emptyset\}$,
and the elements $u_1,\ldots,u_k$ are constructed on the previous step.

From $\ker \psi\subseteq \ker \varphi$,
we have the injective embedding $U_{\RB\As}(C)$ into $\RB\As(A)$.

It remains to show, using \eqref{eq:RB-PreOperations},
that the complement $E'$ of the set of all $E$-words
in the base of $\RB\As(A)$ is linearly expressed via $E$-words.
Applying the inductions in $\RB\As(A)$ on the $R$-degree and the degree
of base words, the relations
\begin{gather}
R(a)u = \begin{cases}
(a\succ b)u', & u = bu', b\in B, \label{AsTheo1} \\
R(R(a)t)u' + R(aR(t))u', & u = R(t)u'; \end{cases} \\
aR(bR(u)) = aR(b)R(u) - aR(R(b)u ) = (a\prec b)R(u) - aR(R(b)u),\quad a,b\in B, \label{AsTheo2}
\end{gather}
the relations for $uR(a)$ and $R(R(u)b)a$ analogous to \eqref{AsTheo1}, \eqref{AsTheo2},
and the relations derived from \eqref{eq:corr1}--\eqref{eq:corr4}
by the removing the symbol $*$,
we prove that the elements of $E'$ are linearly expressed via $E$-words.
The theorem is proved.

{\bf Example 5}.
Let $C$ denote a free preassociative algebra generated by a set $X$.
Universal enveloping algebra $U_{\RB}(C)$ is a free associative RB-algebra
generated by the $X$. This statement proven in \cite{FardGuo07}
one could theoretically deduce from Theorem 2;
but the proof of such corollary does not seem to be obvious.

Given a preassociative algebra $C$,
$U_0(C)$ denotes a linear span of all $E$-words of zero $R$-degree in $U_{\RB}(C)$.

{\bf Example 6}.
Let $D$ be a vector space of square matrices from $M_n(\Bbbk)$,
then $D$ is a preassociative algebra with respect
to the operations $a\succ b=ab$ (as in $M_n(\Bbbk)$), $a\prec b=0$.
From $D^2 = D$ we conclude $U_0(D) = D$.

{\bf Example 7}.
Let $K$ be a vector space of the dimension $n^2$ over the field $\Bbbk$
with the operations $\succ,\prec$ defined as $a\succ b = a\prec b = 0$.
We have that $U_0(K)$ equal $\As\langle K\rangle$, a free associative algebra
generated by $K$.

{\bf Corollary 1}.
The pair of varieties $(\RB\As, \pre\As)$ is not a PBW-pair.

{\sc Proof}.
Universal enveloping associative RB-algebras of
finite-dimensional preassocitive algebras $D$ and $K$ (from Examples 6 and 7)
of the same dimension are not isomorphic, else the spaces
$U_0(D)$ and $U_0(K)$ were isomorphic as vector spaces. But we have
$$
\dim U_0(D) = \dim D = n^2<\dim U_0(K) = \dim\As\langle K\rangle = \infty.
$$

Therefore, the structure of universal enveloping associative RB-algebra
of a preassociative algebra $C$ essentially depends on the operations $\succ,\prec$ on $C$.

\section{Universal enveloping Rota---Baxter algebra of post\-associative algebra}

In the paragraph, we will construct universal enveloping associative RB-algebra
for a postassociative algebra $\langle C,\succ,\prec,\cdot\rangle$.
Let $B$ be a linear base of $C$.

Define by induction $E$-words, a subset of the standard base of $\RB_1\As\langle B\rangle$:

1) elements of $B$ are $E$-words of the type 1;

2) given $E$-word $u$, we define $R(u)$ as an $E$-word of the type 2;

3) the word
$$
v = u_0 R(v_1) u_1 R(v_2) u_2 \ldots u_{k-1} R(v_k) u_k,
\quad \deg v\geq 2,\,\deg_R(v)\geq1,
$$
is $E$-word of the type 3, if
$u_1,\ldots, u_{k-1}\in B$, $u_0,u_k\in B\cup\emptyset$,
$v_1,\ldots,v_{k}$ are $E$-words of any type, moreover,
$R(v_2),\ldots,R(v_{k-1})$ are semigood elements,
$R(v_1)$ is right-good.
Given $u_0\neq\emptyset$, $R(v_1)$ is semigood.
Given $u_k\neq\emptyset$, $R(v_k)$ is semigood,
else $R(v_k)$ is left-good.

{\bf Theorem 3}.
The set of all $E$-words forms a linear base of
universal enveloping associative RB-algebra of $C$.

{\bf Lemma 6}.
Let $D$ denote a linear span of all $E$-words.
One can define such bilinear operation $*$ on the space $D$ that
($k$--$l$ denotes below the condition on the product $v*u$,
where $v$ is an $E$-word of the type $k$ and $u$ is an $E$-word of the type $l$.)

1--1: given $v,u\in B$, we have
\begin{equation}\label{tr-def:1--1}
v*u = v\cdot u.
\end{equation}

1--2: given $v = a\in B$, an $E$-word $u = R(p)$ of the type 2, we have
\begin{equation}\label{tr-def:1--2}
v*u = \begin{cases}
   a\prec b, & p = b \in B, \\
  (a\prec b)R(x)-a*R(R(b)*x)-a*R(b*x), & p = bR(x), b\in B, \\
   aR(p), & R(p)\ \mbox{is left-good}.
               \end{cases}
\end{equation}

1--3: given $v = a\in B$, an $E$-word $u$ of the type 3, we have
\begin{multline}\label{tr-def:1--3}
v*u \\
= \begin{cases}
  (a\cdot b)R(x)u', & u = bR(x)u',\,b\in B, \\
   (a\prec b)*(R(y)*u')- (a*R(R(b)*y))*u' & u = R(bR(y))u',\,b\in B,\\
   \hfill -(a*R(b*y))*u', & R(bR(y))\ \mbox{is not semigood}, \\
   aR(x)u', & u = R(x)u',\,R(x)\ \mbox{is semigood}.
      \end{cases}
\end{multline}

3--1: given $u = a\in B$, an $E$-word $v$ of the type 3, we have
\begin{equation}\label{tr-def:3--1}
v*u  = \begin{cases}
  v'R(x)(b\cdot a), & v = v'R(x)b,\,b\in B, \\
  (v'*R(y))*(b\succ a)- v'*(R(y*R(b))*a) & v = v'R(R(y)b),\,b\in B, \\
  \hfill  - v'*(R(y*b)*a), & \\
  v'R(x)a, & v = v'R(x),\,R(x)\ \mbox{is good}.
      \end{cases}
\end{equation}

2--2: given $E$-words $v = R(v')$, $u = R(u')$ of the type 2, we have
\begin{equation}\label{tr-def:2--2}
v*u = R(v') * R(u') = R( R(v')*u' + v'*R(u') + v'*u').
\end{equation}

2--3: given an $E$-word $v = R(x)$ of the type 2, an $E$-word $u$ of the type 3, we have
\begin{multline}\label{tr-def:2--3}
v*u \\
 = \begin{cases}
  (a\succ b)R(y)u', & x=a,u = b R(y)u',\,a,b\in B, \\
   R(z)(a\succ b)R(y)u'-R(z*R(a))*u & x=R(z)a,u = b R(y)u',a,b\in B, \\
   \hfill -R(z*a)*u, &  \\
    R(x)u, & u = bu',\,b\in B,\,R(x)\ \mbox{is right-good}, \\
   R( R(x)*y + x*R(y)+x*y)*u', & u = R(y)u'.
      \end{cases}
\end{multline}

3--3: given $E$-words $v,u$ of the type 3, we have
\begin{equation}\label{tr-def:3--3}
v*u=\begin{cases}
  (v'*R(R(x)*y+ x*R(y)))*u'  & v = v' R(x),\,u = R(y)u', \\
  \hfill + (v'*R(x*y))*u', & \\
  v'R(x)(a*u), & v = v'R(x)a,\,u = R(y)u',a\in B, \\
  (v*a)R(y)u', & v = v'R(x),\,u = aR(y)u',a\in B, \\
  v'R(x)(a\cdot b)R(y)u', & v = v'R(x)a,\,u = bR(y)u',\,a,b\in B.
      \end{cases}
\end{equation}

The following conditions are also satisfied:

L1) Given left-good $E$-word $R(x)$, $b\in B$,
and $R(b)*x = \sum u_i+\sum u_i'$, where $u_i,u_i'$ are $E$-words,
$\deg_R(u_i')<\deg_R(u_i) = \deg_R(x)+1$,
we have that $R(u_i)$ is left-good for every $i$.

L2) Given right-good $E$-word $R(y)$,
any $E$-word $x$, and $R(x)*y = \sum u_i+\sum u_i'$, where
$u_i,u_i'$ are $E$-words, $\deg_R(u_i')<\deg_R(u_i) = \deg_R(x)+\deg_R(y)+1$,
we have that $R(u_i)$ is right-good for every $i$.

L3) Given $E$-word $x\not\in B$, an $E$-word $u$ of the type 1 or 3,
and $R(x)*u = \sum u_i+\sum u_i'$, where $u_i,u_i'$ are $E$-words,
$\deg_R(u_i')<\deg_R(u_i) = \deg_R(x)+\deg_R(u)+1$,
we have that $u_i$ is an $E$-word of the type 3
and begins with an $R$-letter, i.e., has a view of $u_i = R(x_i)u_i'$, for every $i$.

L4) Given an $E$-word $u = u'a$, $a\in B$, of the type 1 or 3,
any $E$-word $v$, and $v*u = \sum u_i$, where $u_i$ are $E$-words,
we have that $u_i = u_i'a_i$, $a_i\in B$.

L5) Given $a\in B$, an $E$-word $u = R(x)bu'$,
where the word $bu'$ could be empty, $b\in B$,
we have $a*u = \sum\limits_{j\in J_1} a_j R(x_j)b_ju' + \sum\limits_{j\in J_2} (a_j\cdot b_j)u'$ with
$a_j,b_j\in B$, $j\in J_1\cup J_2$, of the view
$a_j = (a,c_1,c_2,\ldots,c_k)_\prec$,
$b_j = (d_1,d_2,\ldots d_k,b)_\succ$, $c_k,d_k\in B$.
Moreover, the number of summands and values of $c_j,d_k,x_j$
depend only on $x$.

The relations for the products of $E$-words
of types 2--1, 3--2 and the conditions R1--R5
are defined analogously to the products of $E$-words
of types 1--2, 2--3 and the conditions
L1--L5 by the inversion of letters, multipliers
and operations symbols, wherein
$\succ$ and $\prec$ turn in each other,
and $\cdot$ does not change.

{\sc Proof}.
Let us define the operation $*$ with the prescribed conditions
for $E$-words $v,u$ by induction on $r = \deg_R(v)+\deg_R(u)$.
For $r = 0$, define $v*u = v\cdot u$, $v,u\in B$,
it satisfies the condition 1--1.
It is easy to see that all conditions L1--L5 (R1--R5) hold.

The case $r = 1$ is possible only if
$v = a\in B$, $u = R(b)$, $b\in B$ (or $v = R(b)$, $u=a$).
Define $a*R(b) = a\prec b$, $R(b)*a = b\succ a$,
this satisfies the condition 1--2.
It is clear that all conditions L1--L5 (R1--R5)
for $r = 1$ are fulfilled.

Suppose that the product $v*u$ is yet defined for
all pairs of $E$-words $v,u$ such that
$\deg_R(v)+\deg_R(u)<r$, $r\geq2$, and all conditions
of the statement on $*$ are satisfied.
Let us define $*$ on $E$-words $v,u$ with
$\deg_R(v)+\deg_R(u) = r$
by induction on $d = \deg(v)+\deg(u)$.

For $d = 2$, consider the cases 1--2 and 2--2
(the case 2--1 is analogous to 1--2).

1--2: let $v = a\in B$, $u = R(p)$ be an $E$-word of the type,
$\deg_R(u)=r\geq 2$, define $v*u$ by \eqref{tr-def:1--2}.

2--2: let $v = R(v')$, $u = R(u')$ be an $E$-word of the type,
define $v*u$ by \eqref{tr-def:2--2}.

The products $R(b)*x$, $b*x$, $R(v')*u'$, $v'*R(u')$, $v'*u'$ in
\eqref{tr-def:1--2}, \eqref{tr-def:2--2}
are defined by induction on $r$.
The multiplication of $a$ on $R(b*x)$
is defined, as $\deg_R(a)+\deg_R(R(b*x))<r$.
The multiplication of $a$ on $R(R(b)*x)$
is defined by the condition L1
holding by inductive assumption for $R(b)*x$.
The element $p = bR(x)$ is an $E$-word of the type 3,
so $R(x)$ is left-good and, therefore,
the concatenation $(a\prec b)R(x)$ is correct.
The conditions L2--L5 (and R2--R5)
in the cases $d=2$, 2--1 and 2--2 hold;
the conditions L1 and R1
are realizable only in the case 2--2,
hence, they are also fulfilled.

For $d > 2$, define the product $v*u$ for
$E$-words pairs of the cases 1--2, 2--3, 3--3, 1--3, 3--1;
the products for the cases 2--1 and 3--2 are defined
analogously up to the inversion.

1--2: let $v = a\in B$, $u = R(p)$ be an $E$-word of the type 2,
$\deg_R(u)=r\geq 2$, define $v*u$ by \eqref{tr-def:1--2}.

The definition is correct by the same reasons as in the case 1--2 for $d=2$.

2--3: let $v = R(x)$ be an $E$-word of the type 2,
$u$ be an $E$-word of the type 3, define $v*u$ by \eqref{tr-def:2--3}.

By the conditions R1, L2, R3 and induction on $r$,
the definition of the product in the case 2--3 is correct.
The definition of $v*u$ is consistent with L1--L5 (R1--R5).

3--3: let $v,u$ be $E$-words of the type 3,
define $v*u$ by \eqref{tr-def:3--3}.

The correctness of products in \eqref{tr-def:3--3}
follows in the first case by the conditions L2, L3, R2, R3,
in the second and third --- by L4, R4.
The definition of $v*u$ is consistent with L1--L5 (R1--R5).

1--3: let $v=a\in B$, $u$ be an $E$-word of the type 3,
define $v*u$ by \eqref{tr-def:1--3}.

3--1: let $u=a\in B$, $v$ be an $E$-word of the type 3,
define $v*u$ by \eqref{tr-def:3--1}.

The conditions L1 and R1 provide that
the definition in the cases 1--3 and 3--1 are true.
The definition of the product in the case 1--3 for $u = R(bR(y))u'$,
where $R(bR(y))$ is not semigood, is reduced to the cases 3--1 and 3--3.
In the last variant of the case 3--1,
the product is expressed by the one from the case 1--3.
The process of computation $v*u$ is finite by the same reasons
as in the proof of Lemma 1.

Actually the product in the case $k$--$l$ from Lemma 6
differs from the one from Lemma 1 in additional summands of less $R$-degree.

The definition of $v*u$ in the cases 1--3 and 3--1 is consistent with L1--L5 (R1--R5).

{\bf Lemma 7}.
The space $D$ with the operations $*$, $R$ is an RB-algebra.

{\sc Proof}.
It follows from \eqref{tr-def:2--2}.

{\bf Lemma 8}.
The relations
$R(a)*b = a\succ b$, $a*R(b) = a\prec b$, $a*b = a\cdot b$ hold in $D$ for every $a,b\in B$.

{\sc Proof}.
It follows from Lemma 1, equality 1--1 \eqref{tr-def:1--1},
the first case of \eqref{tr-def:1--2}, and analogous relation of 2--1.

{\bf Lemma 9}.
The operation $*$ on $D$ is associative.

{\sc Proof}.
Let us prove associativity on $D$ by inductions on two parameters:
at first, on summary $R$-degree $r$ of the $E$-words triple $x,y,z$,
at second, on summary degree $d$ of the triple $x,y,z$.

For $r = 0$, we have $(x,y,z) = (x\cdot y)\cdot z - x\cdot (y\cdot z) = 0$, $x,y,z\in E_0$,
as the product $\cdot$  is associative in the algebra $A$.

Let $r>0$ and suppose that associativity
for all triples of $E$-words with the less summary $R$-degree is proven.

We prove the statement for the triples $x,y,z$,
in which $y$ is an $E$-word of the type 1 or 3, $d$ is any.
We consider the cases 1--1--2, 1--1--3, and 1--3--1,
all others could be proven analogously.

1--1--2. a) Given $x = a$, $y = b$, $z = R(c)$, $a,b,c\in B$,
by \eqref{id:PostAs}, we compute
$$
(x*y)*z - x*(y*z)
 = (a\cdot b)\prec c - a\cdot(b\prec c) = 0.
$$

b) If $z = R(x)$ is left-good, associativity is obvious.

c) Given $x = a$, $y = b$, $z = R(cR(t))$, $a,b,c\in B$, we have
\begin{multline}\label{PostAs1-1-2c-L}
(x*y)*z = (a\cdot b)*R(cR(t)) \\
= ((a\cdot b)\prec c)R(t) - (a\cdot b)*R(R(c)*t) - (a\cdot b)*R(c*t);
\end{multline}
\begin{multline}\label{PostAs1-1-2c-R}
x*(y*z) = a*((b\prec c)R(t) - b*R(R(c)*t) - b*R(c*t)) \\
= (a\cdot (b\prec c))R(t) - a*(b*R(R(c)*t)) - a*(b*R(c*t))).
\end{multline}

The first summands of RHS of \eqref{PostAs1-1-2c-L}
and \eqref{PostAs1-1-2c-R} equal by \eqref{id:PostAs},
the second ones --- by the condition L1 and induction on $r$,
the third ones --- by induction on $r$.

1--1--3. a) Given $x = a$, $y = b$, $z = cz'$, $a,b,c\in B$,
by \eqref{id:PostAs}, we compute
$$
(x*y)*z - x*(y*z)
 = ((a\cdot b)\cdot c)z' - (a\cdot(b\cdot c))z' = 0.
$$

b)  If $z = R(x)z'$ and $R(x)$ is semigood, associativity is obvious.

c) Given $x = a$, $y = b$, $z = z_0 z'$, $a,b,c\in B$,
where an $R$-letter $z_0$ is not semigood. Let us present $z_0$ as
$$
z_0 = R(a_1 R( R(a_3 R(R(\ldots R(a_{2k-1} R(R(a_{2k+1}R(t))a_{2k}) )  \ldots )a_4) )a_2 ) ),
$$
$R(t)$ is good, $a_1,\ldots,a_{2k+1}\in B$.

Let $k=0$, $z_0 = R(cR(t))$, $R(t)$ be good, we have
\begin{multline}\label{PostAs1-1-3c-L}
(x*y)*z = (a\cdot b)*R(cR(t))z' \\
= ((a\cdot b)\prec c)R(t)*z' - ((a\cdot b)*R(R(c)*t))*z' - ((a\cdot b)*R(c*t))*z';
\end{multline}
\begin{multline}\label{PostAs1-1-3c-R}
x*(y*z) = a*((b\prec c)R(t)*z' - (b*R(R(c)*t))*z' - (b*R(c*t))*z' \\
= (a\cdot (b\prec c))R(t)*z' - a*((b*R(R(c)*t))*z') - a*((b*R(c*t))*z').
\end{multline}

The first summands of RHS of \eqref{PostAs1-1-3c-L}
and \eqref{PostAs1-1-3c-R} equal by \eqref{id:PostAs} and induction on $r$,
the third ones --- by induction on $r$,
the second ones --- by the condition L1, induction on $r$,
and the fact that $(b*R(R(c)*t))*z' = bR(s_0))z' + (b*R(s_1))*z'$
for $\deg_R(s_1)<\deg_R(s_0) = \deg_R(t)+1$.

Let $k>0$. Presenting $z_0 = R(cR(R(t)d))$, by \eqref{PostAs1-1-3c-L},
\eqref{PostAs1-1-3c-R}, and the induction on $r$, we have
\begin{multline}\label{PostAs1-1-3c-As}
(x,y,z)
 = - ((a\cdot b)*R(R(c)*R(t)d))*z' + a*((b*R(R(c)*R(t)d))*z') \\
 = - (a\cdot b)*(R( R(R(c)*t)*d + R(c*R(t))*d + R(c*t)*d  ))*z' \\
 + a*((b*R( R(R(c)*t)*d + R(c*R(t))*d + R(c*t)*d  ))*z' )
\end{multline}

As $z_0$ is not semigood, $t = eR(p)$, $e\in B$,
from \eqref{PostAs1-1-3c-As} we have
$$
(x,y,z)
 = (a\cdot b)*(R( R(c*R(R(e)*p))*d )*z')
 - a*( (b* R( R(c*R(R(e)*p))*d  ))*z' )
$$
Continuing on the process, we obtain
\begin{multline}\label{PostAs1-1-3c-As2}
(x,y,z)
 = \\
[(a\cdot b)*(
R(R(a_1*R( R(a_3*R(\ldots R(a_{2k-1}*R(R(a_{2k+1}*R(q))*a_{2k}))\ldots ))*a_4))*a_2)*z') \\
- a*( (b*R(R(a_1*R( R(a_3*R(\ldots R(a_{2k-1}*R(R(a_{2k+1}*R(q))*a_{2k}))\ldots ))*a_4))*a_2))*z')]
\end{multline}
As $R(q)$ is good, the inner product $a_{2k+1}*R(q)$
equals $s_0+s_1$ with $\deg_R(s_1)<\deg_R(s_0) = \deg_R(q)+1$,
$R(s_0)$ is good. Hence, $(x,y,z) = 0$ by induction on $r$
and the condition 1--3.

1--3--1. The only case which is not analogous to the cases considered above
is the following: $x = a$, $y = R(s)bR(t)$, $z = c$, $a,b,c\in B$,
$R(s)$ is not semigood, $R(t)$ is not right-good.
From one hand,
\begin{multline}\label{PostAs1-3-1L}
(a*R(s)bR(t))*c
 = aR(s_0)bR(t)*c + \sum\limits_i a_i R(s_i)b_iR(t)*c + \sum\limits_{k,i'} (\bar{a}_k\cdot b_{i'})R(t)*c \\
 = aR(s_0)bR(t_0)c + \sum\limits_j aR(s_0)b_j R(t_j)c_j + \sum\limits_{l,j'} aR(s_0)(b_{j'}\cdot \bar{c}_l) \\
 + \sum\limits_i a_i R(s_i)b_iR(t_0)c + \sum\limits_{i,j} a_i R(s_i)(b_i)_j R(t_j)c_j
 + \sum\limits_{i,l,j'} a_i R(s_i)((b_i)_{j'}\cdot \bar{c}_l) \\
 + \sum\limits_{k,i'} (\bar{a}_k\cdot b_{i'})R(t_0)c + \sum\limits_{j,k,i'} (\bar{a}_k\cdot b_{i'})_j R(t_j)c_j
 + \sum\limits_{k,l,i',j'} ((\bar{a}_k\cdot b_{i'})_{j'}\cdot \bar{c}_l),
\end{multline}
where $\deg_R(s_i)<\deg_R(s_0) = \deg_R(s)$, $\deg_R(t_j)<\deg_R(t_0) = \deg_R(t)$,
all variations of letters of $a,b,c$ lie in $B$,
and indexes $i,j,k,l,i',j'$ run over some disjoint finite sets.

From another hand,
\begin{multline}\label{PostAs1-3-1R}
a*(R(s)bR(t)*c)
 = a*R(s)bR(t_0)c + a*\sum\limits_j R(s)b_jR(t_j)c_j + a*\sum\limits_{l,j'} R(s)(b_{j'}\cdot \bar{c}_l) \\
 = aR(s_0)bR(t_0)c + \sum\limits_i a_i R(s_i)b_i R(t_0)c + \sum\limits_{k,i'} (\bar{a}_k\cdot b_{i'})R(t_0)c \\
 + \sum\limits_j a R(s_0)b_jR(t_j)c_j + \sum\limits_{i,j} a_i R(s_i)(b_j)_i R(t_j)c_j
 + \sum\limits_{j,k,i'} (\bar{a}_k\cdot (b_j)_{i'})R(t_j)c_j \\
 + \sum\limits_{l,j'} aR(s_0)(b_{j'}\cdot\bar{c}_l) + \sum\limits_{i,l,j'} a_iR(s_i)(b_{j'}\cdot\bar{c}_l)_i
 + \sum\limits_{k,l,i',j'} (\bar{a}_k\cdot(b_{j'}\cdot \bar{c}_l)_{i'}).
\end{multline}
Comparing \eqref{PostAs1-3-1L} and \eqref{PostAs1-3-1R},
it is enough to state that
\begin{gather}
\sum\limits_{i,j} a_i R(s_i)(b_i)_j R(t_j)c_j
 = \sum\limits_{i,j} a_i R(s_i)(b_j)_i R(t_j)c_j,
\label{PostAs1-3-1As1} \\
\sum\limits_{i,l,j'} a_i R(s_i)((b_i)_{j'}\cdot \bar{c}_l)
 = \sum\limits_{i,l,j'} a_iR(s_i)(b_{j'}\cdot\bar{c}_l)_i,
\label{PostAs1-3-1As2} \\
\sum\limits_{j,k,i'} (\bar{a}_k\cdot b_{i'})_j R(t_j)c_j
 = \sum\limits_{j,k,i'} (\bar{a}_k\cdot (b_j)_{i'})R(t_j)c_j,
\label{PostAs1-3-1As3} \\
\sum\limits_{k,l,i',j'} ((\bar{a}_k\cdot b_{i'})_{j'}\cdot \bar{c}_l)
 = \sum\limits_{k,l,i',j'} (\bar{a}_k\cdot(b_{j'}\cdot \bar{c}_l)_{i'}).
\label{PostAs1-3-1As4}
\end{gather}

The equality \eqref{PostAs1-3-1As1} follows from the conditions L5, R5 and Statement.
The equality \eqref{PostAs1-3-1As2} due to the conditions L5, R5 is equivalent to the equality
$$
(d_1,\ldots,d_p,(b,r_1,\ldots,r_q)_\prec)\cdot \bar{c}_l)_\succ
 = ((d_1,\ldots,d_p,b)_\succ,r_1,\ldots,r_q))_\prec\cdot\bar{c}_l,
$$
which is true by Statement and \eqref{id:PostAs}.
The proof of \eqref{PostAs1-3-1As3} and \eqref{PostAs1-3-1As4} is analogous.

Hence, we have to prove associativity only for triples
$x,y,z$, in which $y$ is an $E$-word of the type 2.
The definition of $*$ by Lemma 6 is symmetric
with respect to the inversion, except the cases 1--3 and 3--1, so
associativity in the triple of types $k$--2--$l$
leads to associativity in the triple $l$--2--$k$.

Prove associativity for $d = 3$ and $E$-word $y$ of the type~2.
We consider only those cases whose proof is not directly
analogous to the one from Lemma 5.

1--2--1. Given $x = a$, $y = R(b)$, $z = c$, $a,b,c\in B$,
by \eqref{id:PostAs}, we compute
$$
(a*R(b))*c - a*(R(b)*c) = (a\prec b)\cdot c - a\cdot (b\succ c) = 0.
$$

1--2--2. a) Let $x = a$, $y = R(b)$, $z = R(c)$, $a,b,c\in B$. By \eqref{id:PostAs}, we have
$$
(a*R(b))*R(c) - a*(R(b)*R(c)) = (a\prec b)\prec c
- a\prec(b\succ c + b\prec c + b\cdot c) = 0.
$$

b) $x = a$, $y = R(b)$, $z = R(cR(t))$, $a,b,c\in B$.
\begin{multline}\label{PostAs1-2-2a-L}
(x*y)*z = (a\prec b)*R(cR(t)) \\
= ((a\prec b)\prec c)R(t) - (a\prec b)*R(R(c)*t) - (a\prec b)*R(c*t).
\end{multline}
\begin{multline}\label{PostAs1-2-2a-R}
x*(y*z) = a*R( R(b)*cR(t) + b*R(cR(t)) + b*cR(t) ) \\
= a*R((b\succ c)R(t)) + a*R(b*R(cR(t))) + a*R((b\cdot c)R(t)) \\
= (a\prec (b\succ c))R(t) - a*R(R(b\succ c)*t)-a*R((b\succ c)*t) \\
 + a*R( (b\prec c)R(t) - b*R(R(c)*t) - b*R(c*t) ) \\
 + (a\prec (b\cdot c))R(t) - a*R(R(b\cdot c)*t)-a*R((b\cdot c)*t).
\end{multline}
Applying inductive assumption,
write down the penultimate string of \eqref{PostAs1-2-2a-R}:
\begin{multline}\label{PostAs1-2-2a-R2}
a*R( (b\prec c)R(t) - b*R(R(c)*t) - b*R(c*t) ) \\
 = (a\prec (b\prec c))R(t) - a*R(R(b\prec c)*t) - a*R((b\prec c)*t) \\
 - (a\prec b)*R(R(c)*t) + a*R(R(b)*R(c)*t) + a*R(b*R(c)*t)  \\
 - (a\prec b)*R(c*t) + a*R(R(b)*c*t) + a*R(b*c*t).
\end{multline}
Substituting \eqref{PostAs1-2-2a-R2} in \eqref{PostAs1-2-2a-R}
and subtracting the result from \eqref{PostAs1-2-2a-L},
by \eqref{id:PostAs}, we have $a*(R(b),R(c),t)$
equal to zero by induction.

c) $x = a$, $y = R(bR(t))$, $z = R(u)$, $a,b\in B$.
Notice that $(a*R(R(b)*t))*R(u) = a*(R(R(b)*t)*R(u))$.
Indeed, let $R(b)*t = s_1+s_2$,
where $s_1$ denotes a linear combination of $E$-words
starting with a $R$-letter and $s_2$ --- starting with a letter from $B$.
Thus, $\deg_R(s_2)<\deg_R(t)+1$ and, hence, by induction, we have
\begin{multline*}
(a*R(R(b)*t))*R(u)
 = (aR(s_1)+a*R(s_2))*R(u) \\
 = a*(R(s_1)*R(u))+a*(R(s_2)*R(u))
 = a*(R(s_1+s_2)*R(u)) \\
 = a*(R(R(b)*t)*R(u)).
\end{multline*}
Applying the result, compute
\begin{multline}\label{PostAs1-2-2b-L}
(x*y)*z = (a*R(bR(t))*R(u)
= ((a\prec b)R(t) - a*R(R(b)*t) - a*R(b*t))*R(u) \\
= (a\prec b)*R(R(t)*u + t*R(u) + t*u) \\
 - a*R( R(b)*t*R(u) + R(R(b)*t)*u + R(b)*t*u ) \\
 - a*R( R(b*t)*u + b*t*R(u) + b*t*u ).
\end{multline}
\begin{multline}\label{PostAs1-2-2b-R}
x*(y*z) = a*R( R(bR(t))*u + bR(t)*R(u) + bR(t)*u ) \\
= a*R( R(bR(t))*u + bR(t)*u ) \\
 + (a\prec b)*R(R(t)*u + t*R(u) + t*u) \\
 - a*R(R(b)*(R(t)*u + t*R(u) + t*u)) \\
 - a*R(b*(R(t)*u + t*R(u) + t*u)).
\end{multline}
Subtracting \eqref{PostAs1-2-2b-R} from \eqref{PostAs1-2-2b-L}
and using the equality $(R(b),R(t),u) = 0$ holding by induction,
we get zero.

For $d>3$, one consider other cases analogously to the cases from the proof of Lemma~5.
Thus, Lemma 9 is proven.

{\sc Proof of Theorem 3}.
Let us prove that the algebra $D$ is
exactly universal enveloping algebra for the preassociative algebra $C$,
i.e., is isomorphic to the algebra
$$
U_{\RB_1\As}(C)
 = \RB_1\As\langle B|a\succ b = R(a)b,\,a\prec b = aR(b),\,a\cdot b = ab,\,a,b\in B\rangle.
$$
By the construction, the algebra $D$ is generated by $B$.
Therefore, $D$ is a homomorphic image of a homomorphism $\varphi$ from $U_{\RB_1\As}(C)$.
We will prove that all basic elements of $U_{\RB_1\As}(C)$ are linearly expressed by $D$,
then $\mathrm{ker}\,\varphi = (0)$ and $D\cong U_{\RB_1\As}(C)$.

In $U_{\RB_1\As}(C)$, the equality $ab = a\cdot b$ holds,
therefore, $U_{\RB_1\As}(C)$ is a subspace of $\RB_1\As\langle B\rangle$.
Denote by $E'$ the complement of the set of all $E$-words
in the base of $\RB_1\As\langle B\rangle$.
Applying the inductions in $\RB_1\As(A)$ on the $R$-degree and the degree
of base words, the equalities $xy = x\cdot y$, $x,y\in B$, the relations
\begin{gather}
R(a)u = \begin{cases}
(a\succ b)u', & u = bu', b\in B, \label{tr-AsTheo1} \\
R(R(a)t)u' + R(aR(t))u' + R(at)u', & u = R(t)u'; \end{cases} \\
aR(bR(u)) = (a\prec b)R(u) - aR(R(b)u ) - aR(bu),\quad a,b\in B, \label{tr-AsTheo2}
\end{gather}
the relations for $uR(a)$ and $R(R(u)b)a$ analogous to \eqref{tr-AsTheo1}, \eqref{tr-AsTheo2},
and the relations derived from the analogues of \eqref{eq:corr1}--\eqref{eq:corr4}
by the removing the symbol $*$, we prove that the elements of $E'$
are linearly expressed via $E$-words.

{\bf Corollary 2}.
The pair of varieties $(\RB_\lambda\As, \post\As)$ is a PBW-pair.

\medskip
The author expresses his gratitude to P. Kolesnikov
for important corrections.

The research is supported by RSF (project N 14-21-00065).

\noindent
Gubarev Vsevolod \\
Sobolev Institute of Mathematics of the SB RAS \\
Acad. Koptyug ave., 4 \\
Novosibirsk State University  \\
Pirogova str., 2 \\
Novosibirsk, Russia, 630090\\
{\it E-mail: wsewolod89@gmail.com}

\end{document}